# Nonparametric goodness-of fit testing in quantum homodyne tomography with noisy data

**Katia Meziani**[*]

*Laboratoire de Probabilités et Modèles Aléatoires, Université Paris VII (Denis Diderot), 75251 Paris Cedex 05, France*
*e-mail:* meziani@math.jussieu.fr

**Abstract:** In the framework of quantum optics, we study the problem of goodness-of-fit testing in a severely ill-posed inverse problem. A novel testing procedure is introduced and its rates of convergence are investigated under various smoothness assumptions. The procedure is derived from a projection-type estimator, where the projection is done in $\mathbb{L}_2$ distance on some suitably chosen *pattern* functions. The proposed methodology is illustrated with simulated data sets.

**AMS 2000 subject classifications:** Primary 62G05, 62G10, 62G20; secondary 81V80.
**Keywords and phrases:** Density matrix, Goodness-of fit test, Minimax rates, Nonparametric test, Pattern Functions estimation, Quantum homodyne tomography, Wigner function.

In quantum mechanics, the results of physical measurements performed on a physical system prepared in a certain state are random. This fact provides a large variety of problems and mathematical and applied statistics offer many methods to answer them. For example, statisticians are asked to decide which measurement to perform in order to collect most informative observations or to rebuild the subjacent quantum object from the resulting observations. In particular, physicists have created laser in view of many applications (medical, ...). An important question is to know whether they have generated the desired one. In order to answer this question, we propose a testing procedure based on the indirect observations they gathered by a measurement called quantum homodyne tomography (QHT). To describe our framework and to motivate our work, we present in Section 1 the needed physical background. In Section 2, we introduce our statistical model and define the nonparametric class containing the density matrices we are dealing with. To make clear to physicists the statistical part, we recall some basic notions on the statistical tests in Section 2.2. The testing procedure and the theoretical results assessing its asymptotic behavior are presented in Section 3. Section 4 contains an illustration of the proposed methodology on simulated data sets while the proof are postponed to Section 5.

---

[*]The authors acknowledge the support of the French Agence Nationale de la Recherche (ANR), under grant StatQuant (JC07 07205763) "Quantum Statistics".

## 1. Physical background

Section 1.1 is a short introduction to quantum optics and we introduce in Section 1.2 the very useful and meaningful *pattern* functions. In Section 1.3, we present some examples of quantum states actually created in laboratory, on which our testing procedure may be applied.

### *1.1. Short introduction to quantum optics*

In quantum mechanics, the quantum state fully describes all aspects of a physical system. Mathematically, the quantum state is described via **density operators** $\rho$ on a complex Hilbert space $\mathcal{H}$ such that

1. $\rho$ is self-adjoint (or Hermitian): $\rho = \rho^*$, where $\rho^*$ is the adjoint of $\rho$,
2. $\rho$ is positive: $\rho \geq 0$, or equivalently $\langle \psi, \rho\psi \rangle \geq 0$ for all $\psi \in \mathcal{H}$,
3. $\rho$ has trace one: $\text{Tr}(\rho) = 1$.

Moreover, to each measurable physical property or quantity corresponds a self-adjoint operator, say $\mathbf{X}$, on the space of states $\mathcal{H}$. This operator is called an *observable*. Unlike in classical mechanics, when performing a certain measurement of an observable $\mathbf{X}$ of a physical system prepared in a quantum state $\rho$, the result is in general random and is described by a probability distribution of a random variable $X$.

In this paper, the studied quantum system is a monochromatic light in a cavity described by the state of a quantum harmonic oscillator. Hereafter, we consider the complex Hilbert space $\mathcal{H} = \mathbb{L}_2(\mathbb{R})$, the space of square integrable complex-valued functions on the real line. In quantum optics, the quantum harmonic oscillator can be represented by two operators on $\mathcal{H} = \mathbb{L}^2(\mathbb{R})$: the position operator $\mathbf{Q}$ (or electric field) and the quantity of movement operator $\mathbf{P}$ (or magnetic field)

$$(\mathbf{Q}\psi_1)(x) = x\psi_1(x), \tag{1}$$

$$(\mathbf{P}\psi_2)(x) = -i\frac{\partial \psi_2(x)}{\partial x}, \tag{2}$$

with $\psi_1$, $\psi_2$ functions of $\mathcal{H}$. The position $\mathbf{Q}$ and momentum $\mathbf{P}$ do not commute and satisfy the canonical commutation relation:

$$[\mathbf{Q}, \mathbf{P}] = \mathbf{Q}\mathbf{P} - \mathbf{P}\mathbf{Q} = i\mathbf{1},$$

where $\mathbf{1}$ stands for the identity operator on $\mathcal{H}$. This relation may be understood as follows. The measurement of the position $\mathbf{Q}$ necessarily disturbs the particle's momentum $\mathbf{P}$, and vice versa. They cannot be measured simultaneously. Thus, one cannot obtain a couple of random variables $(Q, P)$. However, a linear combination

$$\mathbf{X}_\phi = \cos(\phi)\mathbf{Q} + \sin(\phi)\mathbf{P}, \text{ for every phases } \phi \in [0, \pi],$$

can theoretically be measured by a technique put in practice for the first time in [29] and called quantum homodyne tomography (QHT[1]), where the phase $\Phi$ is randomly and uniformly distributed on $[0, \pi]$.

In the *case of noiseless* measurement and for a phase $\Phi = \phi$, the result would be a random variable $X$ with probability density noted $p_\rho(\cdot|\phi)$. Moreover, the characteristic function of $X$ is

$$\mathcal{F}_1[p_\rho(\cdot|\phi)](t) = \mathrm{Tr}\left(\rho e^{it\mathbf{X}_\phi}\right) = \mathrm{Tr}\left(\rho e^{it(\cos(\phi)\mathbf{Q}+\sin(\phi)\mathbf{P})}\right).$$

Furthermore, to each state $\rho$ corresponds a **Wigner function** $W_\rho$, which gives an equivalent representation of the quantum state $\rho$. The associated Wigner function $W_\rho$ can be defined rigorously by its Fourier transform $\mathcal{F}_2[W_\rho]$ with respect to both variables

$$\widetilde{W_\rho}(u,v) := \mathcal{F}_2[W_\rho](u,v) = \mathrm{Tr}\left(\rho \exp(iu\mathbf{Q} + iv\mathbf{P})\right).$$

The Wigner function is a mapping from $\mathbb{R}^2$ to $\mathbb{R}$ such that $\iint W_\rho(q,p)dqdp = 1$. Note that the Wigner function can take negative values. For this reason, $W_\rho$ is a real-valued function regarded as a generalized joint probability density (quasi-probability density) of the two random variables $Q$ and $P$ that we would get if we could measure simultaneously the two observables $\mathbf{Q}$ and $\mathbf{P}$, which are respectively given by equations (1) and (2). It is well-known that its Radon transform $\Re[W_\rho]$ is such that

$$p_\rho(x|\phi) = \Re[W_\rho](x,\phi). \tag{3}$$

Here, $\Re$ denotes the Radon transform, taking functions $W_\rho(q,p)$ on $\mathbb{R}^2$ into functions $\Re[W_\rho](x,\phi)$ on $\mathbb{R} \times [0,\pi]$ formed by integration along lines of direction $\phi$ and distance $x$ from the origin, expressed as

$$\Re[W_\rho](x,\phi) = \int_{-\infty}^{\infty} W_\rho(x\cos\phi - t\sin\phi,\ x\sin\phi + t\cos\phi)dt.$$

As $\Phi$ is chosen uniformly and independently on $[0,\pi]$, we can define the probability density function $p_\rho(x,\phi)$ of $(X,\Phi)$ w.r.t. the measure $\frac{1}{\pi}\lambda$, where $\lambda$ stands for the Lebesgue measure on $\mathbb{R} \times [0,\pi]$ by

$$p_\rho(x,\phi) = \frac{1}{\pi}\mathcal{R}[W_\rho](x,\phi). \tag{4}$$

### *1.2. Fock basis and pattern functions*

An important equation in physics, especially in quantum mechanics, is the Schrödinger equation. The time-independent Schrödinger equation is written as

$$\left[-\frac{1}{2}\frac{\partial^2}{\partial x^2} + \frac{x^2}{2}\right]\psi = w\psi, \quad w \in \mathbb{R}, \tag{5}$$

---

[1] We refer the interested reader to the book [18] or the paper [16] for further details on QHT.

where $w$ is the energy level. It turns out that the energies are "quantized", and may only take discrete values noted $w_k = k + 1/2$, $k \in \mathbb{N}$. For a given frequency $w_k$, there are two fundamental solutions: $\psi_k$ and $\varphi_k$. One is a normalized function, the function $\psi_k$, which is such that $\int \psi_k^2 = 1$ and is called the regular wave function. The other one, the function $\varphi_k$, is called the irregular one as it cannot be normalizable as $\psi_k$ is.

The functions $\{\psi_k\}_{k \in \mathbb{N}}$ form a orthonormal basis of $\mathbb{L}_2(\mathbb{R})$. This particular basis, physically very meaningful, is called the **Fock** basis and is written as follows:

$$\psi_k(x) := \frac{1}{\sqrt{\sqrt{\pi}2^k k!}} H_k(x) e^{-x^2/2}. \tag{6}$$

Here, $H_k(x) := (-1)^k e^{x^2} \frac{d^k}{dx^k} e^{-x^2}$ denotes the Hermite polynomial of degree $k$.

The density operator defined in the previous section corresponds to a **density matrix** under some orthonormal basis. In the Fock basis, the matrix entries $\rho_{j,k}$ of the state $\rho$ can be expressed as expected values of functions $F_{j,k}(X_\ell, \Phi_\ell) = f_{j,k}(X_\ell) e^{-i(k-j)\Phi_\ell}$, where $f_{j,k} = f_{k,j}$ are bounded real functions called *pattern functions* [18]. For all $j, k \in \mathbb{N}$,

$$\rho_{j,k} = \iint_0^\pi p_\rho(x, \phi) f_{j,k}(x) e^{-i(k-j)\phi} d\phi dx,$$

where $p_\rho(x, \phi)$ is the joint probability density of $(X, \Phi)$. In other words

$$\rho_{j,k} = E_\rho [F_{j,k}(X, \Phi)]. \tag{7}$$

Equation (7) expresses the idea that one can reconstruct any density matrix element $\rho_{j,k}$ using the *pattern* functions. First introduced in [20], the *pattern* functions $f_{j,k}$ are well known in physics and are defined in [19] as the first derivatives of products of the two fundamental solutions $\psi_k$ and $\varphi_k$ of the Schrödinger equation given in (5) for $j \geq k$

$$f_{j,k}(x) = \frac{\partial}{\partial x} \left(\psi_j(x) \varphi_k(x)\right).$$

A concrete expression for their Fourier transform using generalized Laguerre polynomials can be found in [28, 4]. For $j \geq k$

$$\tilde{f}_{j,k}(t) = \pi(-i)^{j-k} \sqrt{\frac{2^{k-j} k!}{j!}} |t| t^{j-k} e^{-\frac{t^2}{4}} L_k^{j-k}(\frac{t^2}{2}), \tag{8}$$

where $\tilde{f}_{j,k}$ denotes the Fourier transform of the *pattern* function $f_{j,k}$ and $L_k^\alpha(x)$ denotes the generalized Laguerre polynomial. We note that the *pattern* functions $f_{j,k}(x)$ are even functions for even differences $j - k$ and odd functions for odd ones

$$f_{j,k}(-x) = (-1)^{j-k} f_{j,k}(x).$$

TABLE 1
*Examples of quantum states*

**Vacuum state**
- $\rho_{0,0} = 1$ rest zero,
- $p_\rho(x|\phi) = e^{-x^2}/\sqrt{\pi}$.

**Single photon state**
- $\rho_{1,1} = 1$ rest zero,
- $p_\rho(x|\phi) = 2x^2 e^{-x^2}/\sqrt{\pi}$.

**Coherent-$q_0$ state** $q_0 \in \mathbb{R}$
- $\rho_{j,k} = e^{-|q_0/\sqrt{2}|^2}(q_0/\sqrt{2})^{j+k}/\sqrt{j!k!}$,
- $p_\rho(x|\phi) = \exp(-(x - q_0 \cos(\phi))^2)/\sqrt{\pi}$.

**Squeezed state** $(M, \delta) \in \mathbb{R}_+^2$, $\xi \in \mathbb{R}$
- $\rho_{j,k} = C(M, \xi)(\frac{1}{2}\tanh(\xi))^{k+j} H_j(\delta) H_k(\delta)/\sqrt{j!k!}$,
- $p_\rho(x|\phi) = \exp\left[\left(\sin^2(\phi)(xe^{-2\xi}\cos(\phi) - (x\cos(\phi) - \alpha)e^{2\xi})^2\right)/\left(e^{2\xi}\sin^2(\phi) + e^{-2\xi}\cos^2(\phi)\right)\right.$
$\left. \times -e^{-2\xi}(x\cos(\phi) - \alpha)^2 - e^{2\xi}x^2\sin^2(\phi)\right]/\sqrt{\pi(e^{2\xi}\sin^2(\phi) + e^{-2\xi}\cos^2(\phi))}$.

**Thermal state** $\beta > 0$
- $\rho_{j,k} = (1 - e^{-\beta})e^{-\beta k}\mathbb{1}_{j=k}$,
- $p_\rho(x|\phi) = \sqrt{\tanh(\beta/2)/\pi}\exp(-x^2\tanh(\beta/2))$.

**Schrödinger cat-$q_0$** $q_0 > 0$
- $\rho_{j,k} = 2(q_0/\sqrt{2})^{j+k}/\left(\sqrt{j!k!}(\exp(q_0^2/2) + \exp(-q_0^2/2))\right)$, for $j$ and $k$ even, rest zero,
- $p_\rho(x|\phi) = \left(\exp(-(x - q_0\cos(\phi))^2) + \exp(-(x + q_0\cos(\phi))^2)\right.$
$\left. +2\cos(2q_0 x\sin(\phi))\exp(-x^2 - q_0^2\cos^2(\phi))\right)/\left(2\sqrt{\pi}(1 + \exp(-q_0^2))\right)$.

## *1.3. Examples of quantum states*

We present in Table 1 examples of pure quantum states, which can be created at this moment in laboratory and belong to the class $\mathcal{R}(B, r, L)$ with $r = 2$. A state is called pure if it cannot be represented as a mixture (convex combination) of other states, i.e., if it is an extreme point of the convex set of states. This is equivalent to the density matrix being a one dimensional projector, i.e., of the form $\rho = \mathbf{P}_\psi$ for some unit vector $\psi$. Equivalently, a state $\rho$ is pure if $\text{Tr}(\rho^2) = 1$. All other states are called mixed states.

Let us discuss these few examples of quantum states. Among the pure states we consider the *single photon* state and the *vacuum* state, which is the pure state of zero photons. Note that the *vacuum* state would provide a random variable of Gaussian probability density $p_\rho(x|\phi)$ via the ideal measurement of quantum homodyne tomography. We consider also the *coherent-$q_0$* state, which characterizes the laser pulse with the number of photons Poisson distributed with an average of $M$ photons. The *Squeezed* states (see e.g. [7]) have Gaussian Wigner functions whose variances in the two directions have a fixed product. The parameters $M$ and $\xi$ are such that $M \geq \sinh^2(\xi)$, $C(M, \xi)$ is a normalization constant, $\alpha = ((M - \sinh^2(\xi))^{1/2})/(\cosh(\xi) - \sinh(\xi))$, and $\delta = (\alpha/(\sinh(2\xi)))^{1/2}$. The *Schrödinger cat* state is described by a linear superposition of two *coherent* vectors (see e.g. [25]). Table 1 gives some explicit density matrix coefficients $\rho_{j,k}$ and probability densities $p_\rho(x|\phi)$.

## 2. Problem formulation

In this paper, the studied quantum system is a monochromatic light in a cavity, whose state is described by an infinite **density matrix** $\rho$ on the Hilbert space $\mathcal{H} = \mathbb{L}_2(\mathbb{R})$. In this setting, a convenient representation of a quantum state can be obtained by the projection onto the orthonormal Fock basis $(\psi_k)_{k \in \mathbb{N}}$ defined in (6). In quantum optics, physicists produce quantum state of light and via QHT measurements, they gather independent identically distributed random variables containing information on the unknown, underlying quantum state $\rho$. In an ideal framework, the results of measurements would be $(X_\ell, \Phi_\ell)_{\ell=1,\ldots,n}$, independent identically distributed random variables with values in $\mathbb{R} \times [0, \pi]$, with $p_\rho(x, \phi)$ the probability density function of $(X_1, \Phi_1)$ defined in(4).

### 2.1. Statistical model

In this paper we consider a more realistic model in presence of an additional independent Gaussian noise. In practice from $n$ identical, independent **QHT** measurements, we do not collect data $(X_\ell, \Phi_\ell)$, but we observe $(Y_\ell, \Phi_\ell)_{\ell=1,\ldots,n}$ independent identically distributed random variables, as error is add such that

$$Y_\ell := \sqrt{\eta} X_\ell + \sqrt{(1-\eta)/2}\, \xi_\ell, \qquad (9)$$

where $\xi_\ell$ is a sequence of independent identically distributed standard Gaussians, independent of all $(X_\ell, \Phi_\ell)$. The detection efficiency parameter $\eta$, $0 < \eta \leq 1$, is a known parameter and $1 - \eta$ represents the proportion of lost photons due to various losses in the measurement process.

First recall that for any functions $f, g : \mathbb{R} \to \mathbb{R}$, we denote by $f * g$ the convolution product

$$f * g(y) = \int f(y-t) g(t) dt.$$

Due to (9), the density $p_\rho^\eta$ of $(Y_\ell, \Phi_\ell)$ is given by the convolution of the density $p_\rho(\cdot/\sqrt{\eta}, \phi)/\sqrt{\eta}$ with $N^\eta$ the centered Gaussian density of variance $(1-\eta)/2$. In other terms

$$p_\rho^\eta(y, \phi) = \left( \frac{1}{\sqrt{\eta}} p_\rho\left(\frac{\cdot}{\sqrt{\eta}}, \phi\right) * N^\eta \right)(y), \quad \forall y \in \mathbb{R},\ \phi \in [0, \pi].$$

This corresponds to a severely ill-posed inverse problem as the additive noise is super-smooth Gaussian. In the Fourier domain, this relation becomes

$$\mathcal{F}_1[p_\rho^\eta(\cdot, \phi)](t) = \mathcal{F}_1[p_\rho(\cdot, \phi)](t\sqrt{\eta}) \widetilde{N}^\eta(t), \qquad (10)$$

where $\mathcal{F}_1[p_\rho^\eta(\cdot, \phi)]$ denotes the Fourier transform with respect to the first variable and $\widetilde{N}^\eta$ the Fourier transform of the Gaussian density $N^\eta$.

We suppose that the unknown state belongs to a natural class of states from the point of view of applications, the class $\mathcal{R}(B, r, L)$ for $B > 0$ and $r \in ]0, 2]$, defined by

$$\mathcal{R}(B, r, L) := \{\rho \text{ quantum state} : |\rho_{j,k}| \leq L \exp(-B(j+k)^{r/2})\}. \qquad (11)$$

Let us note that all the states described in Table 1 belong to this class with $r = 2$.

## 2.2. Statistical tests

An important problem in quantum optics is to check whether the produced light pulse is in the desired known quantum state $\tau$ or not. The purpose of this paper is to answer this question via goodness-of-fit testing.

More precisely, we consider here the problem of nonparametric goodness-of fit testing from the data $(Y_\ell, \Phi_\ell)$ for $\ell = 1, \ldots, n$; i.e. given $\tau \in \mathcal{R}(B, r, L)$, we define the null hypothesis $H_0$ and the alternative hypothesis $H_1$ as follows:

$$\begin{cases} H_0 : & \rho = \tau, \\ H_1(\mathcal{C}, \varphi_n) : & \rho \in \mathcal{R}(B, r, L) \text{ s.t. that } \|\rho - \tau\|_2 \geq \mathcal{C} \cdot \varphi_n, \end{cases}$$

where $\varphi_n$ is a sequence, which tends to 0 when $n \to \infty$. The physical interpretation of such a test is to check whether the produced light pulse is in a known quantum state $\tau$ ($H_0$ is accepted), or not ($H_1$ is accepted). Here, the distance between the unknown state $\rho$ and the presumed state $\tau$ is measured by the squared-$\mathbb{L}_2$-distance:

$$\|\rho - \tau\|_2^2 = \sum_{j,k \geq 0} |\rho_{j,k} - \tau_{j,k}|^2. \tag{12}$$

In nonparametric statistics, different tools are developed to evaluate the accuracy of a testing procedure. First, let us begin by reminding of some basic definitions. There are two important errors made in a statistical decision process:

- First-type error: the test will reject a correct null hypothesis,
- Second-type error: the test will accept a false null hypothesis.

Given a test procedure $\Omega_n$ such that $\Omega_n = 0$ when we accept $H_0$ and $\Omega_n = 1$ when we reject $H_0$ and decide $H_1$, we denote by $\alpha = P_\tau[\Omega_n = 1]$ and $\beta = \beta(\rho) = P_\rho[\Omega_n = 0]$ the probabilities to make a first-type error and a second-type error respectively under $\tau$ defined in $H_0$ and $\rho$ satisfying $H_1(\mathcal{C}, \varphi_n)$. We would like to control the sum of these two probabilities and we do it as described in Definition 1.

**Definition 1.** *For a given $0 < \lambda < 1$, a test procedure $\Omega_n$ satisfies the upper bound (13) for the testing rate $\varphi_n$ over the smoothness class $\mathcal{R}(B, r, L)$ if there exists a constant $\mathcal{C}^* > 0$ such that for all $\mathcal{C} > \mathcal{C}^*$:*

$$\limsup_{n \to \infty} \left\{ P_\tau[\Omega_n = 1] + \sup_{\rho \in H_1(\mathcal{C}, \varphi_n)} P_\rho[\Omega_n = 0] \right\} \leq \lambda, \tag{13}$$

where $P_\tau$ denote the probability under $\rho = \tau$ defined in $H_0$.

In the statistics, we may also compare the power of several testing procedures having first-type error less than some fixed $\alpha$ and choose the most powerful procedure. The *power* of a statistical test is the probability that the test will reject a false null hypothesis (that it will not make a second-type error) and we denote it by $\Pi$ such that

$$\Pi = 1 - \beta = P_\rho[\Omega_n = 1], \text{ under } \rho \text{ satisfying } H_1(\mathcal{C}, \varphi_n). \tag{14}$$

Remark that, the probability of a second-type error decreases as the *power* increases. In Section 4, to evaluate the performance of our testing procedure $\Omega_n$, we estimate empirically the *power* of our test.

### 2.3. Outline of results

The problem of reconstructing the quantum state of a light beam has been extensively studied in quantum statistics and physical literature. Methods for reconstructing a quantum state are based either on density matrix or on Wigner function estimation. The estimation of the density matrix from averages of data has been considered in the framework of ideal detection [13, 12, 20, 3] as well as in the more general case of an efficiency parameter $\eta$ belonging to the interval $]1/2, 1]$ (cf. [10, 14, 11, 27]). The case $\eta \in ]0, 1]$ has been recently treated in [4]. The latter paper provides also rates of convergence in $\mathbb{L}_2$ loss for an estimator of the Wigner function. The problem of pointwise estimation of the Wigner function has been previously studied in [16] for ideal data and in [9] for noisy data. It should be noted that the results of [16] and most part of the results in [9] are asymptotically minimax not only in the rate, but also in the constant.

In the present work, the goodness-of-fit problem in quantum statistics is considered. There is a large literature on nonparametric testing procedures for the goodness-of-fit of probability distributions. First of all, let us mention the family of test procedures built on certain distances between empirical cumulative distribution functions (c.d.f.), such as the Kolmogorov-Smirnov and the Cramer-Von Mises test statistics, for which extensive results in terms of asymptotic efficiency were established [23]. In order to be more sensitive to low-frequency components or narrow bumps for powerful discrimination, test procedures based on distances between densities were proposed, such as the Bickel-Rosenblatt test [6, 2] for the $\mathbb{L}_2$-distance and the test of [1] for the $\mathbb{L}_1$-distance. Theoretical results on such test statistics naturally stems from their nonparametric function estimation counterparts.

In order to compare nonparametric testing procedures, many approaches were proposed, as reviewed in [23, 17]. A common approach is to analyze the power against sequences of local alternatives $\{f_n\}_{n\geq 1}$ of the form $f_n = f_0 + \varphi_n g$ where $g$ is the direction defining the sequence of local alternatives, and where $\varphi_n \to 0$ as $n$ tends to infinity. Typically, while nonparametric test statistics achieve nontrivial power (i.e. the power of the test is strictly larger than the first type error) against directional alternatives $f_n = f_0 + \varphi_n g$ when $\varphi_n = Cn^{-1/2}$, they achieve nontrivial power against nondirectional alternatives $f_n = f_0 + \varphi_n g_n$ only for $\varphi_n$ slower than $n^{-1/2}$. In other words, achieving nontrivial power uniformly against a large class of alternatives comes at the price of a slower rate than for para-metric testing. The minimax distinguishability framework described in a non-asymptotic setting in [5, 15] and in an asymptotic setting in [17], allows to give precise statements about this phenomenon by characterizing the discrimination rate $\varphi_n$ depending on a smoothness index of the class of alternatives one wish to discriminate from the null hypothesis. Sharp minimax results with pointwise and sup-norm distances have been established in [21] for the regression model and in [26] for the Gaussian white noise model for supersmooth functions. In [8], goodness-of-fit testing in the convolution model have been considered and minimax rates for testing in $\mathbb{L}_2$-norm from indirect observations have been established. The first testing procedure adaptive to smoothness of the alternative function was proposed by [30].

In quantum statistics framework, for the problem of discriminating between two different and fixed states we mention [24] among others. They have establish lower bound for the Bayesian error probability.

The remainder of the article is organized as follows. In Section 3, we present our testing procedure and our theoretical results. We study in Section 4 the numerical performance of our testing procedure. The proofs are deferred to Section 5.

## 3. Testing procedure and results

Our testing method is based on a projection method on pattern functions. We describe first the pattern functions adapted to our setting and derive some useful properties. We define and study the test statistics and give the main results concerning the testing procedure.

### *3.1. Pattern functions*

To take into account the detection losses described by the overall efficiency $\eta$, it is necessary to adapt the *pattern* functions. When $\eta \in ]1/2, 1]$, we denote by $f_{j,k}^\eta$ the suitable functions introduced in [4] and defined by their Fourier transform as follows:
$$\tilde{f}_{j,k}^\eta(t) := \tilde{f}_{j,k}(t) e^{\frac{1-\eta}{4\eta} t^2}. \tag{15}$$
In this paper we develop a procedure for $\eta \in ]1/2, 1]$, but it is possible to get quite similar results with the same procedure for $0 < \eta \leq \frac{1}{2}$ by using modified *pattern* functions, those introduced in [4]. We restrict our study to the more interesting case $\eta \in ]1/2, 1]$ as in practice $\eta$ is around 0.9. The following lemma provides useful upper bounds on the $\mathbb{L}_2$-norm of the $f_{j,k}$ and the $\mathbb{L}_\infty$-norm of the $f_{j,k}^\eta$. From now on, we denote $\gamma := \frac{1-\eta}{4\eta}$.

**Lemma 1** (Aubry *et al.* [4])**.** *There exists a constant $\mathcal{C}_2$ such that for $N$ large enough*
$$\sum_{j,k=0}^{N} \|f_{j,k}\|_2^2 \leq \mathcal{C}_2 N^{\frac{17}{6}}.$$

For $\eta \in ]1/2, 1[$ and $\gamma = \frac{1-\eta}{4\eta}$, there exists constant $\mathcal{C}_\infty^\eta$ such that for $N$ large enough

$$\sum_{j,k=0}^{N} \left\|f_{j,k}^\eta\right\|_\infty^2 \leq \mathcal{C}_\infty^\eta N^{-\frac{2}{3}} e^{16\gamma N}.$$

For the $\mathbb{L}_\infty$-norm of the $f_{j,k}^\eta$, this lemma is slightly different from Lemmata 4 and 5 in [4] where the sum is over $j + k = 0, \ldots, N$. The proof remains similar.

### 3.2. Testing procedure

In this part, we propose a testing procedure which allow to choose among the hypothesis $H_0$ and $H_1$ defined in Section 2.2 by

$$\begin{cases} H_0: & \rho = \tau, \\ H_1(\mathcal{C}, \varphi_n): & \rho \in \mathcal{R}(B, r, L) \text{ s.t. } \|\rho - \tau\|_2 \geq \mathcal{C} \cdot \varphi_n, \end{cases}$$

As in the alternative $H_1$, the true value of $\|\rho - \tau\|_2^2$ is unknown, we have to derive an estimator of this quantity. Then, in a second time, we provide our testing procedure $\Omega_n$.

For the known efficiency parameter $\eta \in ]1/2, 1]$, we define an estimator $M_n$, also called a test statistic of $\|\rho - \tau\|_2^2$. the estimator $M_n$ is as an U-statistic of order 2 based on indirect observations $(Y_i, \Phi_i)_{i=1,\ldots,n}$. For $N := N(n) \to \infty$ as $n \to \infty$, $M_n$ is defined by

$$M_n := \frac{1}{n(n-1)} \sum_{j,k=0}^{N-1} \sum_{\ell \neq m=1}^{n} \left[F_{j,k}^\eta\left(\frac{Y_\ell}{\sqrt{\eta}}, \Phi_\ell\right) - \tau_{j,k}\right] \left[\overline{F_{j,k}^\eta}\left(\frac{Y_m}{\sqrt{\eta}}, \Phi_m\right) - \overline{\tau_{j,k}}\right], \tag{16}$$

where $F_{j,k}^\eta(x, \phi) := f_{j,k}^\eta(x) e^{-i(k-j)\phi}$ uses the *pattern* functions defined in (15) and $\overline{a}$ denote the complex conjugate of $a$. Let us note that the density matrices are infinite, thus we truncate the sum to $N$. The parameter $N$ is called the bandwidth and has to be optimized. We study the properties of $M_n$ in Section 3.3.

Let us discuss the construction of our estimator $M_n$. First remark that we do not estimate

$$d = \sum_{j,k=0}^{\infty} |\rho_{j,k} - \tau_{j,k}|^2 = \sum_{j,k=0}^{N-1} |\rho_{j,k} - \tau_{j,k}|^2 + \sum_{j,k=N}^{\infty} |\rho_{j,k} - \tau_{j,k}|^2 := d_N + R_N$$

but the truncated part $d_N$. For an appropriate choice of $N$, the residual term $R_N$ should be small due to the considered decrease condition (11). Moreover, we have underlined previously that $\rho_{j,k} = E_\rho[F_{j,k}(X, \Phi)]$, thus we can estimate $\rho_{j,k}$ by the empirical version $\hat{\rho}_{j,k}^N$, $0 \leq j, k < N$ (see [4]). For the case $\eta \in ]0, 1[$ we shall use the *pattern* functions defined in (15). We estimate then $d_N$ by an

U-statistic of order 2 given by (16). By the plug-in method, $d_N$ can be estimated by

$$\hat{d}_N = \sum_{j,k=0}^{N-1} |\hat{\rho}_{j,k}^N - \tau_{j,k}|^2 = \sum_{j,k=0}^{N-1} \left| \frac{1}{n} \sum_{\ell=1}^{n} F_{j,k}^\eta \left( \frac{Y_\ell}{\sqrt{\eta}}, \Phi_\ell \right) - \tau_{j,k} \right|^2.$$

However, it is well-known in statistics that such estimators have a too large bias term.

Now, for a constant $\mathcal{C} > 0$ and some threshold $t_n > 0$ defined later, we can define our testing procedure, based on the test statistic $M_n$ defined in (16), as

$$\Omega_n = \mathbb{1}(|M_n| > \mathcal{C} \cdot t_n^2). \tag{17}$$

Here, $\Omega_n = 0$ when $\rho = \tau$ and $\Omega_n = 1$ when the distance between the unknown state $\rho$ and the presumed state $\tau$ is larger than $\mathcal{C} t_n^2$. We study the upper bounds of our procedure $\Omega_n$ in the sense of Definition 13 in Section 3.4.

### 3.3. Properties of the estimator $M_n$

We first remark that each element of the density matrix $\rho_{j,k}$, such that $j, k < N$, is estimated with no bias. Moreover, the estimator $M_n$ is unbiased under $\rho = \tau$ defined in $H_0$ (Remark 1). Remark 2 derives useful tools for the proof of the upper bounds in the sense of Definition 1. From now on, we denote by $E_\rho$ the expected value under $\rho$ satisfying $H_1(\mathcal{C}, \varphi_n)$ and $E_\tau$ under $\rho = \tau$ defined in $H_0$.

**Remark 1.** *Note that for $0 \leq k \leq j \leq N-1$ and $\eta \in ]1/2, 1]$*

$$E_\rho \left[ F_{j,k}^\eta \left( \frac{Y}{\sqrt{\eta}}, \Phi \right) \right] = \rho_{j,k} \quad \text{and} \quad E_\tau \left[ F_{j,k}^\eta \left( \frac{Y}{\sqrt{\eta}}, \Phi \right) \right] = \tau_{j,k}.$$

*Indeed, from Plancherel formula*

$$E_\rho \left[ F_{j,k}^\eta \left( \frac{Y}{\sqrt{\eta}}, \Phi \right) \right] = \iint_0^\pi \int f_{j,k}^\eta \left( \frac{y}{\sqrt{\eta}} \right) e^{-i(k-j)\phi} p_\rho^\eta(y, \phi) d\phi dy$$

$$= \iint_0^\pi f_{j,k}^\eta(y) e^{-i(k-j)\phi} \sqrt{\eta} p_\rho^\eta(\sqrt{\eta} y, \phi) d\phi dy$$

$$= \iint_0^\pi e^{-i(k-j)\phi} \frac{1}{2\pi} \tilde{f}_{j,k}(t) e^{\frac{1-\eta}{4\eta} t^2} \mathcal{F}_1[\sqrt{\eta} p_\rho^\eta(\cdot \sqrt{\eta}, \phi)](t) d\phi dt.$$

*From equation (10), and from Plancherel formula*

$$E_\rho \left[ F_{j,k}^\eta \left( \frac{Y}{\sqrt{\eta}}, \Phi \right) \right] = \iint_0^\pi \frac{e^{-i(k-j)\phi}}{2\pi} \tilde{f}_{j,k}(t) e^{\frac{1-\eta}{4\eta} t^2} \mathcal{F}_1[p_\rho(\cdot, \phi)](t) \widetilde{N}^\eta(t/\sqrt{\eta}) d\phi dt$$

$$= \iint_0^\pi e^{-i(k-j)\phi} f_{j,k}(x) p_\rho(x, \phi) d\phi dx = \rho_{j,k}.$$

*Thus,*

$$E_\rho[M_n] = \sum_{j,k=0}^{N-1} (\rho_{j,k} - \tau_{j,k})^2 \quad \text{and} \quad E_\tau[M_n] = 0.$$

**Remark 2.** *For $N$ large enough and $\rho$ belonging to the class $\mathcal{R}(B, r, L)$*

$$\sum_{j,k=0}^{N-1} E_\rho \left[ \left| F_{j,k}^\eta \left( \frac{Y}{\sqrt{\eta}}, \Phi \right) \right|^2 \right] \leq \begin{cases} \mathcal{C}_\infty^\eta N^{-2/3} e^{16\gamma N} & \eta \in ]1/2, 1[, \\ C \cdot \mathcal{C}_2 N^{17/6} & \eta = 1, \end{cases}$$

where $\mathcal{C}_\infty^\eta$ and $\mathcal{C}_2$ are constants defined in Lemma 1 and $C$ is a positive constant.

*Proof.* Since $E_\rho\bigl[\bigl|F_{j,k}^\eta\bigl(\frac{Y}{\sqrt{\eta}},\Phi\bigr)\bigr|^2\bigr] = E_\rho\bigl[\bigl|f_{j,k}^\eta\bigl(\frac{Y}{\sqrt{\eta}},\Phi\bigr)\bigr|^2\bigr]$ For $\eta \in ]1/2, 1]$

$$\sum_{j,k=0}^{N-1} E_\rho \left[ \left| F_{j,k}^\eta \left( \frac{Y}{\sqrt{\eta}}, \Phi \right) \right|^2 \right] \leq \begin{cases} \sum_{j,k=0}^{N-1} \left\| f_{j,k}^\eta \right\|_\infty^2 & \eta \in ]1/2, 1[, \\ C \sum_{j,k=0}^{N-1} \left\| f_{j,k}^\eta \right\|_2^2 & \eta = 1. \end{cases}$$

For $\eta = 1$, we first apply Lemma 6 in [4] where it has been established that

$$\sup_{x \in \mathbb{R}} \int_0^\pi p_\rho(x, \phi) d\phi \leq C,$$

where $C$ is a positive constant. Hence, the result is a direct consequence of Lemma 1. □

In order to choose the optimal the bandwidth $N$, we do the classical bias/variance trade-off of the risk $E\bigl[\bigl|M_n - \|\rho - \tau\|_2^2\bigr|^2\bigr] = B^2(M_n) + V(M_n)$, where the bias and the variance term are respectively defined as follows:

$$B(M_n) := \bigl| E[M_n] - \|\rho - \tau\|_2^2 \bigr|,$$
$$V(M_n) := E\left[ |M_n - E[M_n]|^2 \right].$$

From now on, we denote by $B_\rho$, $V_\rho$ the variance and the bias terms under $\rho$ satisfying $H_1(\mathcal{C}, \varphi_n)$ and $V_\tau$, $B_\tau$ the variance and the bias terms under $\rho = \tau$ defined in $H_0$. Note that under $\rho = \tau$ defined in $H_0$, the bias term is written $B_\tau(M_n) := |E_\tau[M_n]|$ and is equal to 0. Hence, the following Propositions 1 and 2 evaluate these quantities.

**Proposition 1.** *For $r \in ]0, 2]$, and $\eta \in ]1/2, 1]$ we have*

$$B_\tau(M_n) = 0,$$
$$B_\rho(M_n) \leq C_B N^{2-r/2} e^{-2BN^{r/2}}, \text{ for } N \text{ large enough },$$

*where $C_B$ is a positive constant depending only on $B$, $r$ and $L$.*

**Proposition 2.** *For $r \in ]0, 2]$, $\eta = 1$ and for $N$ large enough such that $N^{17/6}/n \to 0$ as $n \to \infty$, we have*

$$V_\rho(M_n) \leq \frac{8C \cdot \mathcal{C}_2}{n} N^{17/6},$$
$$V_\tau(M_n) \leq \frac{(C \cdot \mathcal{C}_2)^2}{n^2} N^{17/3}.$$

For $r \in ]0,2]$, $\eta \in ]1/2,1[$, $\gamma = (1-\eta)/(4\eta)$ and for $N$ large enough such that $N^{-2/3}e^{16\gamma N}/n \to 0$ as $n \to \infty$, we have

$$V_\rho(M_n) \leq \frac{8\mathcal{C}_\infty^\eta N^{-2/3}}{n} e^{16\gamma N},$$

$$V_\tau(M_n) \leq \frac{(\mathcal{C}_\infty^\eta)^2 N^{-4/3}}{n^2} e^{32\gamma N}.$$

The constants $\mathcal{C}_\infty^\eta$ and $\mathcal{C}_2$ are defined in Lemma 1 and $C$ is positive constant.

### 3.4. Main results

In the following Theorem, we establish upper bounds for the testing rates in the sense of Definition 1 where $\rho$ is the unknown density matrix supposed to belong to the class $\mathcal{R}(B,r,L)$ defined in (11). Theorem 1-(1) deals with the ideal detection case while Theorem 1-(2) and Theorem 2 take into account the Gaussian noise.

**Theorem 1.** *The test procedure $\Omega_n$ defined in (17) for the bandwidth $N(n)$, the threshold $t_n$ and the constant $\mathcal{C}^*$ satisfies the upper bound (13) for the rate $\varphi_n$ such that*

1. *for $r \in ]0,2]$, $\eta = 1$, the bandwidth $N(n) := N_1$ is equal to*

$$N_1 := \left(\frac{\log n}{4B} + \frac{(\log \log n)^2}{4B}\right)^{2/r}, \tag{18}$$

*and the rate*

$$\varphi_n^2 = t_n^2 = n^{-1/2}(\log n)^{\frac{17}{6r}}, \tag{19}$$

2. *for $r = 2$, $\eta \in ]1/2,1[$ and $\gamma := \frac{1-\eta}{4\eta}$, the bandwidth $N(n) := N_2$ is equal to*

$$N_2 := \frac{\log(n)}{4(4\gamma + B)}\left(1 + \frac{8}{3}\frac{\log(\log n)}{\log(n)}\right), \tag{20}$$

*and the rate*

$$\varphi_n^2 = t_n^2 = \log(n)^{\frac{12\gamma - B}{3(4\gamma+B)}} n^{-\frac{B}{2(4\gamma+B)}}. \tag{21}$$

**Theorem 2.** *For $r \in ]0,2[$, $\eta \in ]1/2,1[$ and $\gamma = \frac{1-\eta}{4\eta}$, the test procedure $\Omega_n$ defined in (17) for the bandwidth $N := N_3$ solution of the equation*

$$16\gamma N_3 + 4BN_3^{r/2} = \log n, \tag{22}$$

*the threshold $t_n$ and the constant $\mathcal{C}^*$ satisfies the upper bound (13) for the rate $\varphi_n$*

$$\varphi_n^2 = t_n^2 = N_3^{2-r/2} e^{-2BN_3^{r/2}}. \tag{23}$$

We can remark that

$$\limsup_{n\to\infty}\left\{P_\tau[\Omega_n=1]+\sup_{\rho\in H_1(\mathcal{C},\varphi_n)}P_\rho[\Omega_n=0]\right\}=0.$$

Theorem 1 and 2 provide upper bounds for our testing procedure defined in equation (17) for the testing rates $\varphi_n$. We first remark that our procedure gives nearly parametric rate up to a logarithmic factor when we are in the framework of ideal detection (Theorem 1-(2): no noise) and supersmooth corresponding Wigner functions for all $r \in ]0, 2]$. In the setting of Theorem 1-(2) and Theorem 2, there are good reasons to believe that our testing procedure achieves optimal rates. This remark follows from a recent work of Butucea ([8]), where the author establishes minimax rates for testing in $\mathbb{L}_2$-norm from indirect observations. However, we do not attempt to go that far in this paper.

In the framework of quantum statistics, one can investigate another testing approach based on kernel type estimator of the Wigner function describing the quantum state. Such a testing procedure can been directly derived from the kernel estimator in [22] of the quadratic functional $\iint W_\rho^2$, where $W_\rho$ is the Wigner function associated to the quantum state $\rho$. Our test problem is equivalent to the following:

$$\begin{cases}H_0: & W_\rho=W_{\rho_0},\\ H_1: & \sup_{\rho\in\mathcal{R}(B,r,L)}\|W_\rho-W_{\rho_0}\|_2\geq\mathcal{C}\cdot\varphi_n\end{cases}$$

where $\varphi_n$ is a sequence which tends to 0 when $n \to \infty$. We conjecture that its performances are comparable to those found in this paper and we will leave this analysis for a separate work.

## 4. Simulations

From now on, we set $r = 2$. The purpose of this section is to implement our testing procedure and to investigate its numerical performances. Our motivation is, given a density matrix $\tau \in \mathcal{R}(B, r)$, to decide whether $H_0$ or $H_1$ is accepted

$$\begin{cases}H_0: & \rho=\tau,\\ H_1(\mathcal{C},\varphi_n): & \rho\in\mathcal{R}(B,r,L)\text{ s.t. }\|\rho-\tau\|_2\geq\mathcal{C}\cdot\varphi_n,\end{cases}$$

where $\varphi_n$ is a sequence, which tends to 0 when $n \to \infty$.

We propose to simulate two different situations. In the first one, case **A**, we consider quantum states easily distinguishable, while in the second one, case **B**, we deal with quantum states, which are quite similar and it is difficult to differentiate between them. For $\tau$ defined in $H_0$, we sample our procedure from different density matrices $\rho$ satisfying $H_1(\mathcal{C}, \varphi_n)$ such that in

- the case **A**: $\tau$ is the *vacuum* state, while
  - **a)** $\rho$ is the *vacuum* state ($\rho = \tau$),
  - **b)** $\rho$ is the *single photon* state,
  - **c)** $\rho$ is the *Schrödinger cat*-3 state,

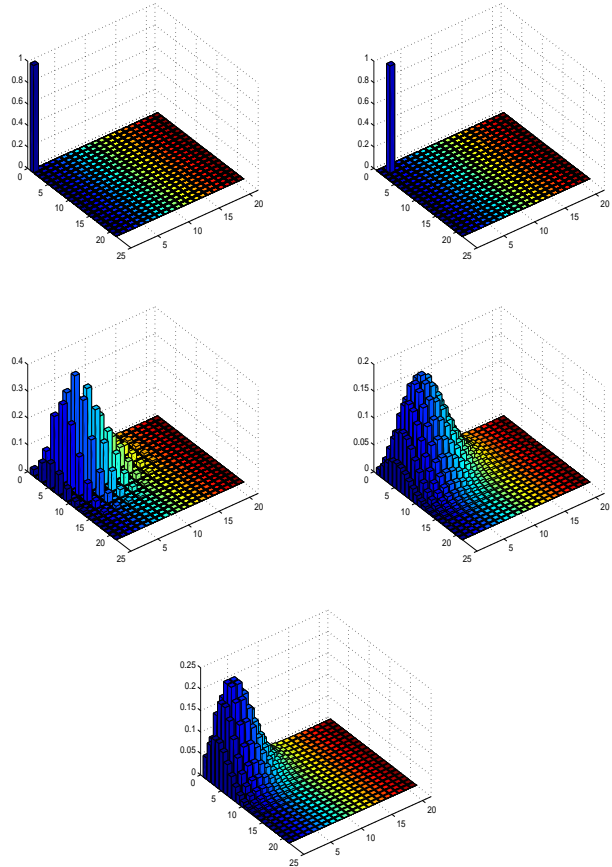

FIGURE 1. *Density matrices of the vacuum state, the single photon state, the Schrödinger cat-3, the coherent-3 state and the coherent-$\sqrt{6}$ state respectively.*

- the case **B**: $\tau$ is the *coherent*-3 state, while
    - **a)** $\rho$ is the *coherent*-3 state ($\rho = \tau$),
    - **b)** $\rho$ is the *coherent*-$\sqrt{6}$ state,
    - **c)** $\rho$ is the *Schrödinger cat*-3 state.

The latter case is a more complicated situation for two main reasons. The *Schrödinger cat-$q_0$* state corresponds to a linear superposition of two *coherent* states ($\pm q_0$). Moreover, the probability density of a *coherent-$q_0$* state is Gaussian with a mean proportional to $q_0$ (see Table 1). Figure 1 represents the density matrices of the states we consider in our simulations. In each situation (cases **A** and **B**), we shall distinguish the case $\eta = 1$, which corresponds to the ideal detection, from the case $\eta = 0.9$, i.e. when noise is present. The latter case is the practical one in laboratory. From now on and for both situations

- when $\eta = 1$, we set $N = 15$,
- when $\eta = 0.9$, we set $N = 14$ and $N = 13$.

The values of $N$ have been chosen after different simulations for different values of $N$. We notice that, when we are in presence of noise ($\eta = 0.9$), we have to take $N$ smaller than the one we choose in the ideal setting. It can be explained by the fact that the variance term of our estimator increase exponentially with $N$ in the case $\eta = 0.9$ (see Proposition 2).

To implement our procedure, we first compute our modified *pattern* function $f_{j,k}^\eta(x)$ in Section 4.1. We implement in Section 4.2 the estimator $M_n$ defined in (16) and finally study the performance of our test procedure $\Omega_n$ defined in (17) in Section 4.3.

## 4.1. Pattern functions $f_{j,k}^\eta$

To implement our procedure, we need the modified *pattern* function $f_{j,k}^\eta(x)$ defined in (15) for all $j \geq k$. For this purpose, we compute here the inverse Fourier Transform $f_{j,k}^\eta(x)$ of the $\tilde{f}_{j,k}^\eta(t)$ given explicitly by (15) and the generalized Laguerre polynomials. Previous authors have used a different method to implement the *pattern* functions in [19, 3] via the following recurrence relation given in [18]:

$$f_{j,k}(x) = 2x\psi_k(x)\varphi_j(x) - \sqrt{2(k+1)}\psi_{k+1}(x)\varphi_j(x) - \sqrt{2(j+1)}\psi_k(x)\varphi_{j+1}(x),$$

for $j \geq k$, otherwise $f_{j,k}(x) = f_{k,j}(x)$. We display in Figure 2-(a) the corresponding graphical representations of the *pattern* functions up to a constant $\pi$ as in physical literature $\pi^{-1}f_{j,k}$ instead of $f_{j,k}$ are often called *pattern* functions. Figure 2-(b) represents some modified *pattern* functions $\pi^{-1}f_{j,k}^\eta$ for $\eta = 0.9$.

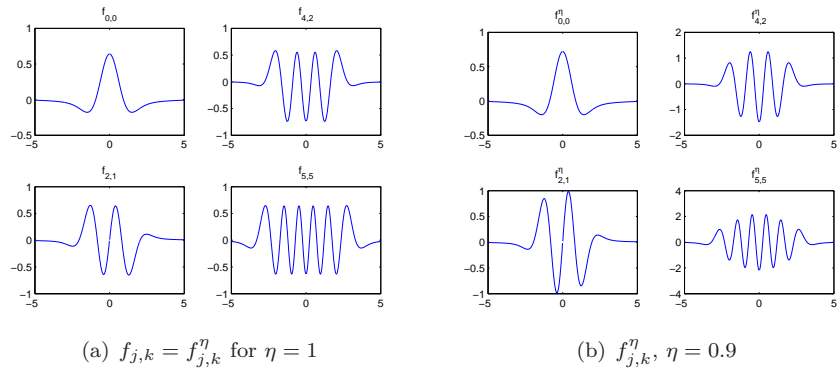

(a) $f_{j,k} = f_{j,k}^\eta$ for $\eta = 1$  (b) $f_{j,k}^\eta$, $\eta = 0.9$

FIGURE 2. *Examples of pattern functions.*

For some of $f_{j,k}$ we expressed below their explicit form.

$$f_{0,0}(x) = 2 - 2e^{-x^2}x\sqrt{\pi}\texttt{Erfi}(x),$$

$$f_{2,1}(x) = e^{-x^2}\left[-2e^{x^2}x(-3+2x^2) + \sqrt{\pi}(1-8x^2+4x^4)\texttt{Erfi}(x)\right],$$

$$f_{4,2}(x) = \frac{e^{-x^2}}{2\sqrt{3}}\left[2e^{x^2}\left(-4+27x^2-24x^4+4x^6\right)\right.$$
$$\left.+\sqrt{\pi}x\left(21-74x^2+52x^4-8x^6\right)\texttt{Erfi}(x)\right],$$

$$f_{5,5}(x) = \frac{e^{-x^2}}{30}\left[2e^{x^2}\left(-30+435x^2-865x^4+526x^6-116x^8+8x^{10}\right)\right.$$
$$\left.+\sqrt{\pi}x\left(225-1425x^2+2160x^4-1160x^6+240x^8-16x^{10}\right)\texttt{Erfi}(x)\right].$$

### 4.2. Implementation of $M_n$

The purpose of this section is to investigate the performance of the estimator $M_n$ of $\|\rho - \tau\|_2^2$ defined in (16) by

$$M_n = \frac{1}{n(n-1)}\sum_{j,k=0}^{N-1}\sum_{\ell\neq m=1}^{n}\left[F_{j,k}^\eta\left(\frac{Y_\ell}{\sqrt{\eta}},\Phi_\ell\right) - \tau_{j,k}\right]\left[\overline{F_{j,k}^\eta}\left(\frac{Y_m}{\sqrt{\eta}},\Phi_m\right) - \overline{\tau_{j,k}}\right],$$

where $F_{j,k}^\eta(x,\phi) = f_{j,k}^\eta(x)e^{-i(k-j)\phi}$ and $\tau$ is defined in $H_0$. Table 2 gives the actual values of $\|\rho - \tau\|_2^2$ in all considered setups.

From now on and for all considered cases (**A** and **B**), the procedure to simulate one value of $M_n$ is designed as follows. We first sample $n = 50\,000$ i.i.d. data from each of the three probability densities $p_\rho$ defined in Table 1 for each considered $\rho$ (Figure 1). The Gaussian noise component $\xi$ is simulated independently with a variance equal to $\sqrt{(1-\eta)/2}$. Then, we run independently of the other runs, 1 000 values of $M_n$ denoted by $\left\{M_n^k\right\}_{k=1,\ldots,1\,000}$, where each value $M_n^k$ is based on $n = 50\,000$ i.i.d. noisy data $(Y_\ell,\Phi_\ell)$.

We recall that the choice of $N = 15$ corresponds to the setting of ideal detection with $\eta = 1$, while $N = 14$ and $N = 13$ deal with the noisy detection with $\eta = 0.9$. Note that for $\rho = \tau$ we expect the $M_n^k$ to be close to 0, while for $\rho$ satisfying $H_1(\mathcal{C},\varphi_n)$ we expect the $M_n^k$ to be close to $\|\rho - \tau\|_2^2 > 0$.

In frameworks **A** and **B**, we draw boxplots of the $\left\{M_n^k\right\}_{k=1,\ldots,1\,000}$ in Figures 3 to 8. For a clearer reading, these boxplots (Figures 3 to 8) have been moved in an appendix at the end of the paper in Section 6. We summarize the median values in Table 3 and we can make the following remarks.

TABLE 2
The values of $\|\rho - \tau\|_2^2$ for $\tau$ the vacuum state (case **A**) and $\tau$ the coherent-3 state (case **B**)

| | **a)** $\rho$ vacuum | **b)** $\rho$ single photon | **c)** $\rho$ Schrödinger cat-3 |
|---|---|---|---|
| Case **A**: | 0 | 2 | 1.9556 |

| | **a)** $\rho$ coherent-3 | **b)** $\rho$ coherent-$\sqrt{6}$ | **c)** $\rho$ Schrödinger cat-3 |
|---|---|---|---|
| Case **B**: | 0 | 0.2812 | 0.9999 |

TABLE 3
The median values of $M_n^k$ for $\tau$ the vacuum state (case **A**) and $\tau$ the coherent-3 state (case **B**)

| Case **A**: | **a)** $\rho$ vacuum | **b)** $\rho$ single photon | **c)** $\rho$ Schrödinger cat-3 |
|---|---|---|---|
| $N = 15$, $\eta = 1$ | 0.0043 | 1.9972 | 1.9417 |
| $N = 14$, $\eta = 0.9$ | 0.0256 | 1.9936 | 1.9360 |
| $N = 13$, $\eta = 0.9$ | 0.0144 | 1.9943 | 1.9378 |
| Case **B**: | **a)** $\rho$ coherent-3 | **b)** $\rho$ coherent-$\sqrt{6}$ | **c)** $\rho$ Schrödinger cat-3 |
| $N = 15$, $\eta = 1$ | $-6.8721 \cdot 10^{-4}$ | 0.2688 | 0.8962 |
| $N = 14$, $\eta = 0.9$ | -0.0466 | 0.1704 | 0.7735 |
| $N = 13$, $\eta = 0.9$ | -0.0376 | 0.1847 | 0.7840 |

TABLE 4
The empirical values of MSE for $\tau$ the vacuum state (case **A**) and $\tau$ the coherent-3 state (case **B**)

| Case **A**: | **a)** $\rho$ vacuum | **b)** $\rho$ single photon | **c)** $\rho$ Schrödinger cat-3 |
|---|---|---|---|
| $N = 15$, $\eta = 1$ | $2.3988 \cdot 10^{-5}$ | $1.3175 \cdot 10^{-4}$ | $4.7877 \cdot 10^{-4}$ |
| $N = 14$, $\eta = 0.9$ | 0.0181 | 0.0056 | 0.0071 |
| $N = 13$, $\eta = 0.9$ | 0.0069 | 0.0025 | 0.0051 |
| Case **B**: | **a)** $\rho$ coherent-3 | **b)** $\rho$ coherent-$\sqrt{6}$ | **c)** $\rho$ Schrödinger cat-3 |
| $N = 15$, $\eta = 1$ | $1.4303 \cdot 10^{-6}$ | $1.9473 \cdot 10^{-4}$ | 0.0110 |
| $N = 14$, $\eta = 0.9$ | 0.0044 | 0.0132 | 0.0531 |
| $N = 13$, $\eta = 0.9$ | 0.0019 | 0.0093 | 0.0484 |

- In the setting **A**: the procedure $M_n$ is a good estimation of $\|\rho - \tau\|_2^2$ both for $\eta = 1$ and $\eta = 0.9$.

- In the setting **B**: when $\rho = \tau$, the procedure $M_n$ shows excellent results since the $M_n^k$ are close to $\|\rho - \tau\|_2^2 = 0$ when either $\eta = 1$ or $\eta = 0.9$. In the ideal case and when $\rho$ is the *coherent*-$\sqrt{6}$ state, the procedure $M_n$ gives a good result since the median of the boxplot of the $M_n^k$ is equal to 0.2688 and the true value $\|\rho - \tau\|_2^2 = 0.2812$. Otherwise, the procedure $M_n$ under evaluates the distance $\|\rho - \tau\|_2^2$.

In order to evaluate the quality of our procedure $M_n$ we estimate the mean square error MSE $= E_\rho[|M_n - \|\rho - \tau\|_2^2|^2]$ as the average over the 1 000 independent runs of $|M_n - \|\rho - \tau\|_2^2|^2$. In other words the MSE is empirically assessed as

$$\frac{1}{1\,000} \sum_{k=1}^{1\,000} \left|M_n^k - \|\rho - \tau\|_2^2\right|^2,$$

with $\|\rho - \tau\|_2^2$ given by Table 2. Table 4 summarizes the results. As we have already noticed the procedure $M_n$ gives excellent results in every cases but it has a larger MSE when we evaluate the distance between the *coherent*-3 state and the *Schrödinger cat*-3 state: MSE is equal to 0.0531 and 0.0484.

### 4.3. Studies of the performance of our test procedure $\Omega_n$

In this part, we would like to confirm the performance of our testing procedure $\Omega_n$. In order to appreciate it, we are interested in the *power* of our test $\Pi$, such that $\Pi = P_\rho[\Omega_n = 1]$ under $\rho$ satisfying $H_1(\mathcal{C}, \varphi_n)$. In this view, we want to compare our estimator $M_n$ with a threshold $\nu_n$ (with $\nu_n = \mathcal{C}^* t_n^2$ in (14)) and decide as follows:

- if $|M_n| > \nu_n$, we accept $H_1$,

- otherwise we accept $H_0$.

In the statistical literature, the parameter $\nu_n = \nu_n(\alpha)$ is also called the critical value associated to $\alpha = P_\tau[\Omega_n = 1]$, which is the probability error of first type (under $\tau$ defined in $H_0$) and defined in Section 2.2. In our framework, we set $\alpha = 1\%$ and $\alpha = 5\%$. As we don't know the density probability of $M_n$ under $\rho = \tau$ defined in $H_0$, we shall evaluate empirically the threshold $\nu_n$ at the testing level $\alpha$ for our particular choices of $\tau$, $n$, $\alpha$ and $\eta$. In this purpose and independently of the future runs, we first compute, for $n = 50\,000$ and for all considered cases 1 000 independent values $\{M_n^k\}_{k=1,\ldots,1\,000}$ of the estimator $M_n$ as described previously. We summarize our results in Table 5. We report that the obtained values $\nu_n$ are larger when $\eta = 0.9$ than when $\eta = 1$. It is due to the noise effect.

From now on, we fix the value of $\nu_n$ as in Table 5. With the same protocol as above and for all considered cases, we compute 1 000 other independent values $\{M_n^k\}_{k=1,\ldots,1\,000}$ of the estimator $M_n$, for $n = 50\,000$. Hence, we evaluate the empirical power of our testing procedure and the empirical first type error. Tables 6 and 7 provide the empirical results obtained by our test procedure in the experiments we deal with. We see that our testing procedure provides very good results even in the framework **B**, we obtain powers $\Pi$ equal to 1 since $N = 13$ when we are in presence of noise. Otherwise, for $N = 14$ and $\eta = 0.9$ the powers of our testing procedure is a little bit degraded, but still remarkably good, since $\Pi = 0.7160$ for $\alpha = 1\%$ and $\Pi = 0.9440$ for $\alpha = 5\%$ in the framework **B** when $\rho$ is the *coherent-$\sqrt{6}$* state, the optimal $N$ when we are in presence of noise is $N = 13$ with powers of test equal to 1.

TABLE 5
*Empirical values of $\nu_n$ for $\tau$ the vacuum state (case **A**) and $\tau$ the coherent-3 state (case **B**)*

|  |  | $N = 15$<br>$\eta = 1$ | $N = 14$<br>$\eta = 0.9$ | $N = 13$<br>$\eta = 0.9$ |
|---|---|---|---|---|
| Case **A** | $\nu_n(1\%)$ | 0.0096 | 0.4605 | 0.2692 |
|  | $\nu_n(5\%)$ | 0.0079 | 0.2739 | 0.1666 |
| Case **B** | $\nu_n(1\%)$ | 0.0028 | 0.1453 | 0.0798 |
|  | $\nu_n(5\%)$ | 0.0022 | 0.1100 | 0.0713 |

TABLE 6
*Empirical values of the first type error $\alpha$ and the power of the test $\Pi$ over 1 000 runs for $\nu_n$ given in Table 5 for $\tau$ the vacuum state (case **A**)*

|  |  | a) $\rho$ vacuum | b) $\rho$ single ph. | c) $\rho$ Schrödinger C.-3 |
|---|---|---|---|---|
|  |  | $\alpha$ | $\Pi$ | $\Pi$ |
| $N=15, \eta=1$ | $\nu_n(1\%)$ | 0.0150 | 1.0000 | 1.0000 |
|  | $\nu_n(5\%)$ | 0.0570 | 1.0000 | 1.0000 |
| $N=14, \eta=0.9$ | $\nu_n(1\%)$ | 0.0070 | 1.0000 | 1.0000 |
|  | $\nu_n(5\%)$ | 0.0550 | 1.0000 | 1.0000 |
| $N=13, \eta=0.9$ | $\nu_n(1\%)$ | 0.0130 | 1.0000 | 1.0000 |
|  | $\nu_n(5\%)$ | 0.0680 | 1.0000 | 1.0000 |

TABLE 7
*Empirical values of the first type error $\alpha$ and the power of the test $\Pi$ over 1 000 runs for $\nu_n$ given in Table 5 for $\tau$ the coherent-3 state (case **B**)*

|  |  | a) $\rho$ coherent-3 | b) $\rho$ coherent-$\sqrt{6}$ | c) $\rho$ Schrödinger C.-3 |
|---|---|---|---|---|
|  |  | $\alpha$ | $\Pi$ | $\Pi$ |
| $N=15, \eta=1$ | $\nu_n(1\%)$ | 0.0090 | 1.0000 | 1.0000 |
|  | $\nu_n(5\%)$ | 0.0530 | 1.0000 | 1.0000 |
| $N=14, \eta=0.9$ | $\nu_n(1\%)$ | 0.0120 | 0.7160 | 1.0000 |
|  | $\nu_n(5\%)$ | 0.0680 | 0.9440 | 1.0000 |
| $N=13, \eta=0.9$ | $\nu_n(1\%)$ | 0.0230 | 1.0000 | 1.0000 |
|  | $\nu_n(5\%)$ | 0.0620 | 1.0000 | 1.0000 |

## 5. Proof

In this section, we give the proofs of the Theorem 1 and 2 derived in Section 3.4 and the proofs of the Proposition 1 and 2 established in Section 3.3.

### 5.1. Proof of the Theorems

In the paragraphs below, we establish the results of the Theorem 1 and 2 derived in Section 3.4. We begin by Theorem 1-(1).

***Proof of Theorem 1-(1).*** Take $r \in ]0,2]$, $\eta=1$, the bandwidth $N_1$ defined in the equation (18) and $\varphi_n^2$ defined in (19), we have by Proposition 2-2 and Proposition 2-1

$$V_\tau(M_n) \leq \frac{(C \cdot C_2)^2}{n^2} N_1^{17/3} \leq C_{V_\tau} \varphi_n^8, \tag{24}$$

$$V_\rho(M_n) \leq \frac{8C \cdot C_2}{n} N_1^{17/6} \leq C_V \varphi_n^4, \tag{25}$$

where $C_V$ and $C_{V_\tau}$ are positive constants depending only on $B$, $r$ and $L$.

Under $\rho = \tau$ defined in $H_0$, from Proposition 1, equations (24) and (19), we bound from above the first type error as follows:

$$P_\tau[\Omega_n = 1] = P_\tau\left[|M_n| \geq \mathcal{C}^* t_n^2\right] \leq \frac{E_\tau[|M_n|^2]}{\mathcal{C}^{*2} t_n^4} = \frac{V_\tau[M_n]}{\mathcal{C}^{*2} t_n^4} \leq \frac{C_{V_\tau} \varphi_n^4}{\mathcal{C}^{*2}} \to 0,$$

as $n \to \infty$.

On the other hand, under $\rho$ satisfying $H_1(\mathcal{C}, \varphi_n)$, we have the second type error bounded as follows:

$$P_\rho[\Omega_n = 0] = P_\rho\left[|M_n| < \mathcal{C}^* t_n^2\right]$$
$$\leq P_\rho\left[|M_n - E_\rho[M_n]| \geq \|\rho - \tau\|^2 - \mathcal{C}^* t_n^2 - B_\rho[M_n]\right].$$

Moreover, under $\rho$ satisfying $H_1(\mathcal{C}, \varphi_n)$, let $\mathcal{C} = \mathcal{C}^*(1+\delta)$, $\delta \in ]0, 1[$, Proposition 1, equations (29) and (23) imply that $B_\rho \leq \frac{\delta}{2}\mathcal{C}^* \varphi_n^2$ for $n$ large enough and then

$$P_\rho[\Omega_n = 0] \leq P_\rho\left[\frac{|M_n - E_\rho[M_n]|}{\sqrt{V_\rho(M_n)}} \geq \frac{\|\rho - \tau\|^2 - \mathcal{C}^* t_n^2 - B_\rho[M_n]}{\sqrt{C_V}\varphi_n^2}\right]$$
$$\leq P_\rho\left[\frac{|M_n - E_\rho[M_n]|}{\sqrt{V_\rho(M_n)}} \geq \frac{\delta \mathcal{C}^*/2}{\sqrt{C_V}}\right]$$
$$\leq \frac{4C_V}{(\delta \mathcal{C}^*)^2} \leq \frac{\lambda}{2},$$

for $\mathcal{C}^*$ large enough. □

**Proof of Theorem 1-(2).** For $r = 2$, $\eta \in ]1/2, 1[$, the bandwidth $N_2$ defined in the equation (20) and $\varphi_n^2$ defined in (23), we know by Proposition 2-2 and Proposition 2-1 that

$$V_\tau(M_n) \leq \frac{(\mathcal{C}_\infty^\eta)^2}{n^2} N_2^{-4/3} e^{32\gamma N_2} \leq C'_{V_\tau} \varphi_n^8, \tag{26}$$

$$V_\rho(M_n) \leq \frac{8\mathcal{C}_\infty^\eta}{n} N_2^{-2/3} e^{16\gamma N_2} \leq C'_V \varphi_n^4, \tag{27}$$

where $C'_V$ and $C'_{V_\tau}$ are positive constants depending only on $B$, $r$, $L$ and $\eta$.

On one hand, under $\rho = \tau$ defined in $H_0$, by Proposition 1, equations (26) and (21), we have the first type error bounded as follows:

$$P_\tau[\Omega_n = 1] = P_\tau\left[|M_n| \geq \mathcal{C}^* t_n^2\right] \leq \frac{E_\tau[|M_n|^2]}{\mathcal{C}^{*2} t_n^4}$$
$$= \frac{V_\tau[M_n]}{\mathcal{C}^{*2} t_n^4} \leq \frac{C'_{V_\tau} \varphi_n^4}{\mathcal{C}^{*2}} \to 0,$$

as $n \to \infty$.

On the other hand, under $\rho$ satisfying $H_1(\mathcal{C}, \varphi_n)$, the second type error is bounded as follows:

$$P_\rho[\Omega_n = 0] = P_\rho\left[|M_n| < \mathcal{C}^* t_n^2\right]$$
$$\leq P_\rho\left[|M_n - E_\rho[M_n]| \geq \|\rho - \tau\|^2 - \mathcal{C}^* t_n^2 - B_\rho[M_n]\right].$$

From Proposition 1, equations (27) and (21), since we are under $\rho$ satisfying $H_1(\mathcal{C}, \varphi_n)$, we get

$$P_\rho[\Omega_n = 0] \leq P_\rho \left[ \frac{|M_n - E_\rho[M_n]|}{\sqrt{V_\rho(M_n)}} \geq \frac{\|\rho - \tau\|^2 - \mathcal{C}^* t_n^2 - B_\rho[M_n]}{\sqrt{C_V'} \varphi_n^2} \right]$$

$$\leq P_\rho \left[ \frac{|M_n - E_\rho[M_n]|}{\sqrt{V_\rho(M_n)}} \geq \frac{\mathcal{C} - \mathcal{C}^* - C_B}{\sqrt{C_V'}} \right]$$

$$\leq \frac{C_V'}{(\mathcal{C} - \mathcal{C}^* - C_B)^2} \leq \frac{\lambda}{2},$$

for $\mathcal{C} > \mathcal{C}^*$ and $\mathcal{C}^*$ large enough. $\square$

**Proof of Theorem 2.** In the case $r \in ]0, 2[$ and $\eta \in ]1/2, 1[$, for the bandwidth $N_3$ solution of the equation (22) and for $\varphi_n^2$ defined in (23), Proposition 2-2 and Proposition 2-1 give that

$$V_\tau(M_n) \leq \frac{(\mathcal{C}_\infty^\eta)^2}{n^2} N_3^{-4/3} e^{32\gamma N_3} \leq C_{V_\tau}'' \varphi_n^8 \cdot N_3^{\frac{2(3r-14)}{3}} \leq C_{V_\tau}'' \varphi_n^8 \quad (28)$$

$$V_\rho(M_n) \leq \frac{8\mathcal{C}_\infty^\eta}{n} N_3^{-2/3} e^{16\gamma N_3} \leq C_V'' \varphi_n^4 \cdot N_3^{\frac{3r-14}{3}} \quad (29)$$

where $C_V''$ and $C_{V_\tau}''$ are positive constants depending only on $B$, $r$, $L$ and $\eta$.

Under $\rho = \tau$ defined in $H_0$, from Proposition 1, equations (28) and (23), it follows that the first type error satisfies

$$P_\tau[\Omega_n = 1] = P_\tau \left[ |M_n| \geq \mathcal{C}^* t_n^2 \right] \leq \frac{E_\tau[|M_n|^2]}{\mathcal{C}^{*2} t_n^4} = \frac{V_\tau[M_n]}{\mathcal{C}^{*2} t_n^4} \leq \frac{C_{V_\tau}'' \varphi_n^4}{\mathcal{C}^{*2}} \to 0,$$

as $n \to \infty$.

Furthermore, under $\rho$ satisfying $H_1(\mathcal{C}, \varphi_n)$, the second type error is such that

$$P_\rho[\Omega_n = 0] = P_\rho \left[ |M_n| < \mathcal{C}^* t_n^2 \right]$$
$$\leq P_\rho \left[ |M_n - E_\rho[M_n]| \geq \|\rho - \tau\|^2 - \mathcal{C}^* t_n^2 - B_\rho[M_n] \right].$$

From Proposition 1, equations (29) and (23), since we are under $\rho$ satisfying $H_1(\mathcal{C}, \varphi_n)$, we deduce that

$$P_\rho[\Omega_n = 0] \leq P_\rho \left[ \frac{|M_n - E_\rho[M_n]|}{\sqrt{V_\rho(M_n)}} \geq \frac{\mathcal{C} - \mathcal{C}^* - C_B}{\sqrt{C_V'' N_3^{\frac{3r-14}{6}}}} \right] \leq \frac{C_V'' N_3^{\frac{3r-14}{3}}}{(\mathcal{C} - \mathcal{C}^* - C_B)^2} \to 0,$$

as $n \to \infty$, for $\mathcal{C} > \mathcal{C}^*$ and $\mathcal{C}^*$ large enough. $\square$

### 5.2. Proof of the Propositions

In the following paragraphs, we give the proofs of Propositions 1 and 2 established in Section 3.3.

**Proof of Proposition 1.** By Remark 1, for $r \in ]0, 2]$

$$B_\rho(M_n) = \left|E_\rho[M_n] - \|\rho - \tau\|_2^2\right| = \left|E_\rho[M_n] - \sum_{j,k \geq 0} |\rho_{j,k} - \tau_{j,k}|^2\right|.$$

It is easy to see that under $\rho = \tau$ defined in $H_0$, we have $B_\tau(M_n) = 0$. Under $\rho$ satisfying $H_1(\mathcal{C}, \varphi_n)$,

$$B_\rho(M_n) = \left|\sum_{j,k=0}^{N-1} |\rho_{j,k} - \tau_{j,k}|^2 - \sum_{j,k \geq 0} |\rho_{j,k} - \tau_{j,k}|^2\right|$$

$$= 2 \sum_{j=N}^{\infty} \sum_{k=0}^{j-1} |\rho_{j,k} - \tau_{j,k}|^2 + \sum_{j=N}^{\infty} |\rho_{j,j} - \tau_{j,j}|^2$$

$$\leq 4 \sum_{j=N}^{\infty} \sum_{k=0}^{j-1} \left(|\rho_{j,k}|^2 + |\tau_{j,k}|^2\right) + 2 \sum_{j=N}^{\infty} \left(|\rho_{j,j}|^2 + |\tau_{j,j}|^2\right)$$

Since $\tau$ and $\rho$ belong to the class $\mathcal{R}(B, r, L)$ defined in (11), it implies

$$B_\rho(M_n) \leq 8L^2 \sum_{j=N}^{\infty} \sum_{k=0}^{j-1} e^{-2B(j+k)^{r/2}} + 4L^2 \sum_{j=N}^{\infty} e^{-2B(j+j)^{r/2}}$$

$$\leq 8L^2 \int_N^{\infty} (j-1) e^{-2Bu^{r/2}} du + 4L^2 \int_N^{\infty} e^{-2B(2u)^{r/2}} du$$

$$\leq C_B N^{2-r/2} e^{-2BN^{r/2}} + C_B' N^{1-r/2} e^{-2B(2N)^{r/2}},$$

where $C_B$ and $C_B'$ denote positive constants depending only on $B$, $r$ and $L$. As $2^{r/2} > 1$ for all $r > 0$,

$$B_\rho(M_n) \leq C_B N^{2-r/2} e^{-2BN^{r/2}} (1 + o(1)),$$

as $N \to \infty$. □

**Proof of Proposition 2.** By centering variables, we have $M_n - E_\rho[M_n] := F_1 + F_2$ where

$$F_1 := \frac{1}{n(n-1)} \sum_{j,k=0}^{N-1} \sum_{\ell \neq m=1}^{n} \left(F_{j,k}^\eta\left(\frac{Y_\ell}{\sqrt{\eta}}, \Phi_\ell\right) - E_\rho\left[F_{j,k}^\eta\left(\frac{Y}{\sqrt{\eta}}, \Phi\right)\right]\right)$$

$$\times \left(\overline{F_{j,k}^\eta}\left(\frac{Y_m}{\sqrt{\eta}}, \Phi_m\right) - E_\rho\left[\overline{F_{j,k}^\eta}\left(\frac{Y}{\sqrt{\eta}}, \Phi\right)\right]\right)$$

$$F_2 := \frac{2}{n} \sum_{j,k=0}^{N-1} \sum_{\ell=1}^{n} \mathcal{R}e\left(\left(F_{j,k}^\eta\left(\frac{Y_\ell}{\sqrt{\eta}}, \Phi_\ell\right) - E_\rho\left[F_{j,k}^\eta\left(\frac{Y}{\sqrt{\eta}}, \Phi\right)\right]\right)\right)$$

$$\times E_\rho \left[ \overline{F_{j,k}^\eta} \left( \frac{Y}{\sqrt{\eta}}, \Phi \right) - \overline{\tau_{j,k}} \right] \right)$$

where $\mathcal{R}e(z)$ stands for the real part of the complex number $z$.

Then, $V_\rho(M_n) := E_\rho[|M_n - E_\rho[|M_n|]|^2] = E_\rho[|F_1|^2] + E_\rho[|F_2|^2]$. First deal with the first term of the previous sum $E_\rho[|F_1|^2]$

$$E_\rho\left[|F_1|^2\right] = \frac{1}{n(n-1)} E_\rho \left[ \left| \sum_{j,k=0}^{N-1} \left( F_{j,k}^\eta\left(\frac{Y_1}{\sqrt{\eta}}, \Phi_1\right) - E_\rho\left[F_{j,k}^\eta\left(\frac{Y}{\sqrt{\eta}}, \Phi\right)\right] \right) \right. \right.$$

$$\left. \left. \left( \overline{F_{j,k}^\eta}\left(\frac{Y_2}{\sqrt{\eta}}, \Phi_2\right) - E_\rho\left[\overline{F_{j,k}^\eta}\left(\frac{Y}{\sqrt{\eta}}, \Phi\right)\right] \right) \right|^2 \right].$$

By the Cauchy-Schwarz inequality on the sum and as $E[|X - E[X]|^2] \leq E[|X|^2]$:

$$E_\rho\left[|F_1|^2\right] \leq \frac{1}{n^2} \left( \sum_{j,k=0}^{N-1} E_\rho\left[\left|F_{j,k}^\eta\left(\frac{Y}{\sqrt{\eta}}, \Phi\right)\right|^2\right] \right)^2.$$

A direct consequence of Remark 2 is

$$E_\rho\left[|F_1|^2\right] \leq \begin{cases} \dfrac{(\mathcal{C}_\infty^\eta)^2}{n^2} N^{-4/3} e^{32\gamma N} & \text{for } \eta \in ]1/2, 1[, \\ \dfrac{(C \cdot \mathcal{C}_2)^2}{n^2} N^{17/3} & \text{for } \eta = 1. \end{cases}$$

By noticing $|\mathcal{R}e(z)| \leq |z|$, the second term of the sum $E_\rho\left[F_2^2\right]$ is such that

$$E_\rho\left[|F_2|^2\right] \leq \frac{4}{n} E_\rho \left[ \left| \sum_{j,k=0}^{N-1} \left( F_{j,k}^\eta\left(\frac{Y_1}{\sqrt{\eta}}, \Phi_1\right) - E_\rho\left[F_{j,k}^\eta\left(\frac{Y}{\sqrt{\eta}}, \Phi\right)\right] \right) \right. \right.$$

$$\left. \left. \times E_\rho\left[F_{j,k}^\eta\left(\frac{Y}{\sqrt{\eta}}, \Phi\right) - \tau_{j,k}\right] \right|^2 \right]$$

By the Cauchy-Schwarz inequality on the sum and as $E[|X - E[X]|^2] \leq E[|X|^2]$:

$$E_\rho\left[|F_2|^2\right] \leq \frac{4}{n} \left( \sum_{j,k=0}^{N-1} E_\rho\left[\left|F_{j,k}^\eta\left(\frac{Y}{\sqrt{\eta}}, \Phi\right)\right|^2\right] \right)$$

$$\times \left( \sum_{j,k=0}^{N-1} \left| E_\rho\left[F_{j,k}^\eta\left(\frac{Y}{\sqrt{\eta}}, \Phi\right)\right] - \tau_{j,k} \right|^2 \right)$$

$$\leq \frac{4}{n} \|\rho - \tau\|_2^2 \left( \sum_{j,k=0}^{N-1} E_\rho\left[\left|F_{j,k}^\eta\left(\frac{Y}{\sqrt{\eta}}, \Phi\right)\right|^2\right] \right)$$

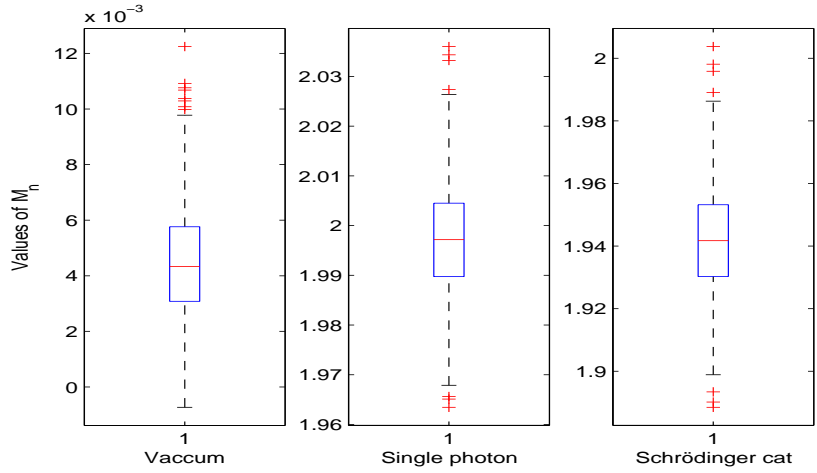

FIGURE 3. *For $\tau$ the vacuum state (case A): **a)** $\rho$ the vacuum state, **b)** $\rho$ the single photon state and **c)** $\rho$ the Schrödinger cat-3 state for $\eta = 1$, $N = 15$.*

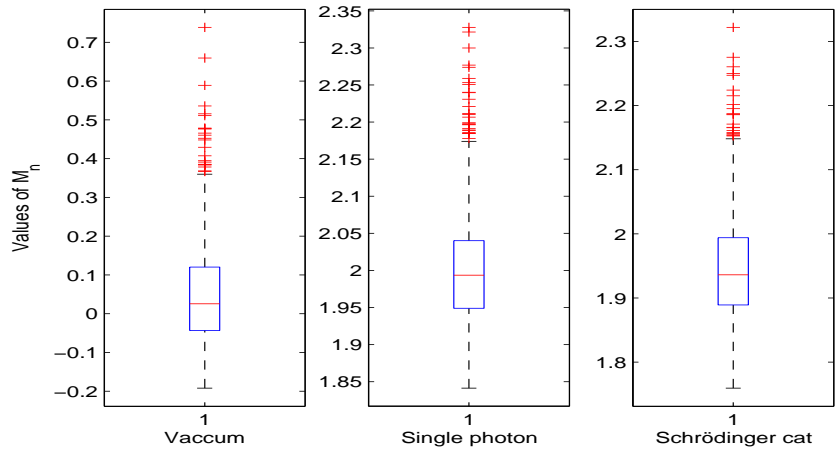

FIGURE 4. *For $\tau$ the vacuum state (case A): **a)** $\rho$ the vacuum state, **b)** $\rho$ the single photon state and **c)** $\rho$ the Schrödinger cat-3 state for $\eta = 0.9$, $N = 14$.*

By Remark 2 and as $\|\rho\|_2^2$, $\|\tau\|_2^2 \leq 1$ we obtain under $\rho$ satisfying $H_1(\mathcal{C}, \varphi_n)$

$$E_\rho\left[|F_2|^2\right] \leq \begin{cases} \dfrac{8\mathcal{C}_\infty^\eta}{n} N^{-2/3} e^{16\gamma N} & \text{for } \eta \in ]1/2, 1[, \\ \dfrac{8C \cdot \mathcal{C}_2}{n} N^{17/6} & \text{for } \eta = 1. \end{cases}$$

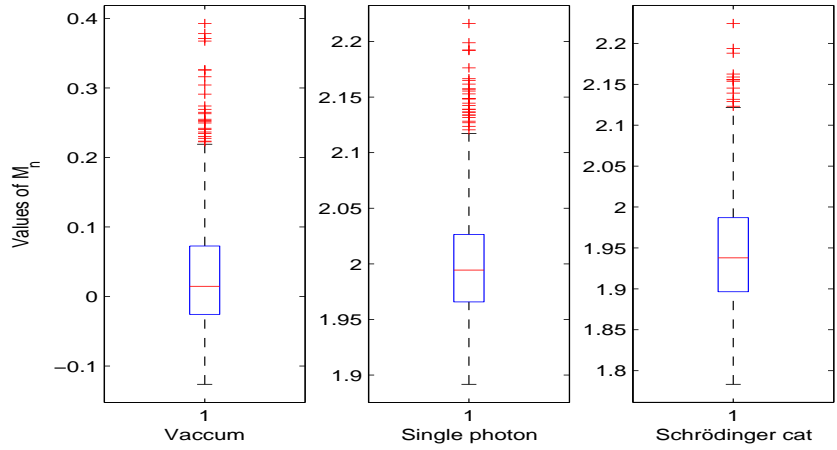

FIGURE 5. *For $\tau$ the vacuum state (case **A**): **a)** $\rho$ the vacuum state, **b)** $\rho$ the single photon state and **c)** $\rho$ the Schrödinger cat-3 state for $\eta = 0.9$, $N = 13$.*

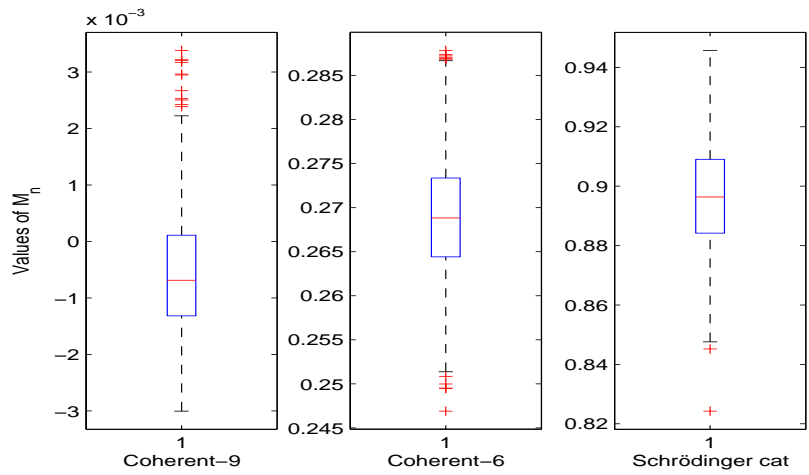

FIGURE 6. *For $\tau$ the coherent-3 state (case **B**): **a)** $\rho$ the coherent-3 state, **b)** $\rho$ the coherent-$\sqrt{6}$ state and **c)** $\rho$ the Schrödinger cat-3 state for $\eta = 1$, $N = 15$.*

If, $N$ is such that the upper bound of $E_\rho\bigl[|F_2|^2\bigr]$ tends to 0 with $n$, then this is the dominant term in the upper bound of $V_\rho(M_n)$. In addition, notice that under $\rho = \tau$ defined in $H_0$, we have $E_\tau\bigl[|F_2|^2\bigr] = 0$. □

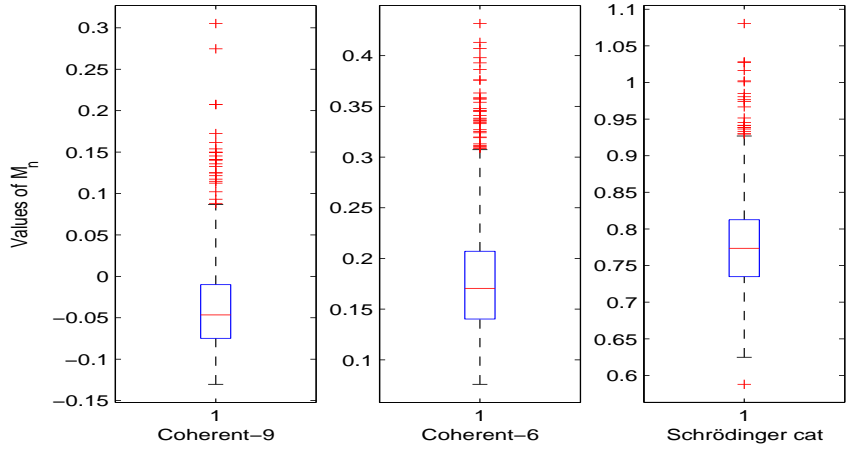

FIGURE 7. *For $\tau$ the coherent-3 state (case **B**): **a)** $\rho$ the coherent-3 state, **b)** $\rho$ the coherent-$\sqrt{6}$ state and **c)** $\rho$ the Schrödinger cat-3 state for $\eta = 0.9$, $N = 14$.*

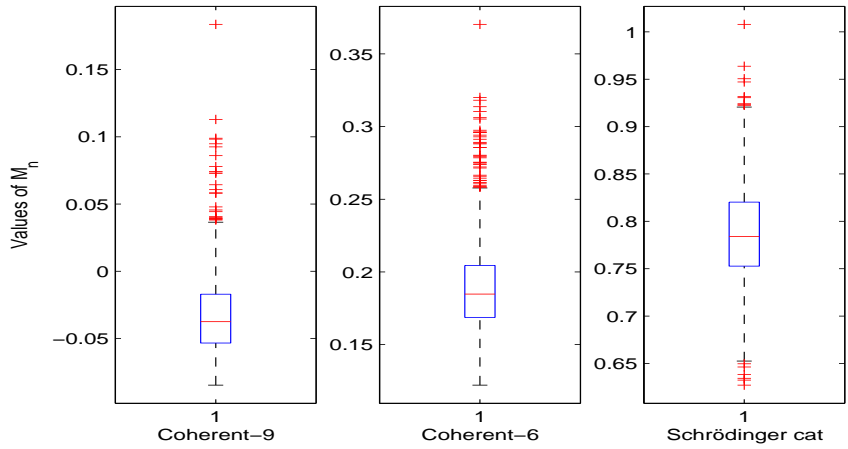

FIGURE 8. *For $\tau$ the coherent-3 state (case **B**): **a)** $\rho$ the coherent-3 state, **b)** $\rho$ the coherent-$\sqrt{6}$ state and **c)** $\rho$ the Schrödinger cat-3 state for $\eta = 0.9$, $N = 13$.*

## 6. Appendix

The boxplots of the 1 000 values $\left\{M_n^k\right\}_{k=1,\ldots,1\,000}$ of the estimator $M_n$ (for $n = 50\,000$) in the case **A** and **B** are represented by Figures 3 to 5 and by Figures 6 to 8 respectively.